\documentclass[11pt]{article}

\makeindex
\usepackage[all]{xy}
\usepackage{graphicx}

\usepackage{amsmath}
\usepackage{amsfonts}
\usepackage{amssymb}
\usepackage{amsthm}
\newcommand{\isom}{\cong}

\newcommand{\Z}{{\bf{Z}}}

\newcommand{\Q}{{\bf{Q}}}

\newcommand{\C}{{\bf{C}}}
\newcommand{\F}{{\bf{F}}}
\newcommand{\T}{{\bf{T}}}

\newcommand{\m}{{\mathfrak{m}}}
\newcommand{\OO}{\mathcal{O}}
\newcommand{\NN}{\mathcal{N}}

\newcommand{\noi}{{\noindent}}

\newcommand{\Mid}{|} %{\mid\!\!}

\newcommand{\miD}{|}%{\!\!\mid}

\newcommand{\divs}{\!\mid\!}
\newcommand{\ndiv}{\!\nmid\!}
\newcommand{\tensor}{\otimes}

\newcommand{\ra}{{\rightarrow}}

 1
\DeclareFontEncoding{OT2}{}{} % to enable usage of cyrillic fonts
  \newcommand{\textcyr}[1]{%
    {\fontencoding{OT2}\fontfamily{wncyr}\fontseries{m}\fontshape{n}%
     \selectfont #1}}
\newcommand{\Sha}{{\mbox{\textcyr{Sh}}}} %{\hbox{\cyr X}}

\newcommand{\ce}{c_{\scriptscriptstyle{E}}}

\newcommand{\caf}{c_{\scriptscriptstyle{A_f}}}

\newcommand{\Agdual}{A_g^{\vee}}
\newcommand{\Af}{{E}}
\newcommand{\Afdual}{E^{\vee}}
\newcommand{\Afp}{E'}
\newcommand{\Ahdual}{A_h^{\vee}}
\newcommand{\Edual}{E^{\vee}}

\newcommand{\Hom}{{\rm Hom}}

\newcommand{\pp}{{\mathfrak{p}}}
\newcommand{\qq}{{\mathfrak{q}}}

\newcommand{\annT}{{{\rm Ann_\T}}}

\newcommand{\tphat}{{\widehat{\T P}}}

\newcommand{\comment}[1]{}
\newcommand{\marginalfootnote}[1]{%
   \footnote{#1}\marginpar{\hfill {\sf\thefootnote}}%
}

\newcommand{\edit}[1]{\marginalfootnote{#1}}

\newtheorem{lem}{Lemma}[section]
\newtheorem{cor}[lem]{Corollary}
\newtheorem{prop}[lem]{Proposition}
\newtheorem{conj}[lem]{Conjecture}
\newtheorem{thm}[lem]{Theorem}

\theoremstyle{definition}

\newtheorem{rmk}[lem]{Remark}

\newcommand{\thetitle}
{A visible factor for analytic rank one}
\begin{document}
%\ssp
\parindent=2em

\title{\thetitle}
\author{Amod Agashe
\footnote{This material is based upon work supported by the National Science 
Foundation under Grant No. 0603668.}}
\maketitle

%\include{abstract}
%abstract.tex
\begin{abstract}
%In an earlier paper, we showed that if vis predicts that q divides Sha,
%it divides Shan. In this paper, we extract a factor of Shan that is
%related to certain congruences, and show that if q divides this factor,
%then q divides Sha, under certain hyp, most serious BSD I.
%Let $N$ be a positive integer, and
%let $f$ be a newform of weight~$2$ on~$\Gamma_0(N)$
%such that its $L$-function~$L(f,s)$ vanishes to order one at~$s=1$.
%Let $\Af$ denote the newform quotient of~$J_0(N)$ associated to~$f$.
Let $E$ be an optimal elliptic curve of conductor~$N$,
such that the $L$-function of~$E$ vanishes to order one at~$s=1$.
Let $K$ be a quadratic imaginary field in which all 
the primes dividing~$N$ split and such that the $L$-function of~$E$ over~$K$
also vanishes to order one at~$s=1$.
%Let $I$ be the index 
%in~$E(K)$ of the subgroup generated by the point in~$E(K)$ obtained from
%the Heegner point on~$X_0(N)(\C)$.
%In view of the Gross-Zagier formula, the second part of
%the Birch and Swinnerton-Dyer conjecture says that up to the Manin constant,
%$I$ is the product of the orders of the arithmetic component groups of~$E$
%and the square root of the order
%of the Shafarevich-Tate group of~$E$ over~$K$.
%The Gross-Zagier theorem gives a formula that expresses
%the Birch and Swinnerton-Dyer conjectural order 
%the Shafarevich-Tate group of~$E$ over~$K$ as a rational number.
In view of the Gross-Zagier theorem, the second part of 
the Birch and Swinnerton-Dyer conjecture says
that the index in~$E(K)$ of the subgroup generated by the Heegner point 
is equal to the product of the Manin constant of~$E$, the Tamagawa numbers
of~$E$, and
the square root of the order of 
the Shafarevich-Tate group of~$E$ (over~$K$).
% in terms of  and the local invariants of~$E$.
% says that up to the Manin constant,
%$I$ is the product of the orders of the arithmetic component groups of~$E$
%and the square root of the order
%of the Shafarevich-Tate group of~$E$ over~$K$.
We extract an integer factor from the index mentioned above and relate
this factor to certain congruences of the newform associated to~$E$ with 
eigenforms of analytic rank bigger than one. We 
use the theory of visibility to show that,
under certain hypotheses (which includes 
the first part of the Birch and Swinnerton-Dyer 
conjecture on rank), if an odd prime~$q$ divides this factor, 
then $q$~divides 
the order of the Shafarevich-Tate group or the order of an
arithmetic component group of~$E$, as predicted by the second
part of the Birch and Swinnerton-Dyer conjecture. 
%as predicted by the second part of the Birch and Swinnerton-Dyer conjecture.
%We extract an integer factor of~$I$ that we relate
%to certain congruences of the newform associated to~$E$ with 
%eigenforms of odd analytic rank bigger than one, and
%use the techniques of visibility to show that,
%under certain hypotheses (which includes 
%the first part of the Birch and Swinnerton-Dyer 
%conjecture on rank), if an odd prime~$q$ divides this factor, then $q^2$~divides 
%the order of the Shafarevich-Tate group,
%as predicted by the second part of the Birch and Swinnerton-Dyer conjecture.
%or the order of a component group of~$A$. 
%that can be related to this product using the theory of visibility,
%assuming the first part of the Birch and Swinnerton-Dyer conjecture.
%We also show that if the newform~$f$ is congruent modulo
%a prime~$p$ to an eigenform of higher analytic rank, then $p$
%divides one of the two factors of~$I$ mentioned above; this 
%is conformity with the fact that in this case,
%by the theory of visibility, $p$ does divide
%the actual order of Shafarevich-Tate group or some component group.
%Thus both of these results 
\end{abstract}

%Warning: This is a very rough draft -- I have mainly just
%jotted down some ideas. At several places,
%one may need minor hypotheses that have been skipped for simplicity.

\section{Introduction and results} \label{sec:intro}

Let $N$ be a positive integer. Let
%, which we might have to assume
%is squarefree in several situations below, 
$X = X_0(N)$ denote the modular curve over~$\Q$
associated to~$\Gamma_0(N)$, and
let $J=J_0(N)$ denote the Jacobian of~$X$, which is
an abelian variety over~$\Q$. 
Let $\T$ denote the Hecke algebra, which is 
the subring of endomorphisms of~$J_0(N)$
generated by the Hecke operators (usually denoted~$T_\ell$
for $\ell \ndiv N$ and $U_p$ for $p\divs N$). 
If $g$ is an eigenform of weight~$2$ on~$\Gamma_0(N)$, then
let $I_g = \annT g$ and let $A_g$ denote the 
quotient abelian variety $J/I_g J$, which is
defined over~$\Q$. 
Also, if $g$ is an eigenform of weight~$2$ on~$\Gamma_0(N)$, then 
the order of vanishing of the $L$-function $L(g,s)$ at~$s=1$ is called
the {\em analytic rank of~$g$}.
Let $f$ be a newform of weight~$2$ on~$\Gamma_0(N)$ 
whose analytic rank is one and which has integer Fourier coefficients.
Then $E = A_f$ is an elliptic curve whose $L$-function vanishes to order
one at $s=1$. 
We denote the quotient map $J \ra J/I_f = \Af$
by~$\pi$. 
%Suppose that $\Af$ has analytic rank one.
%For the sake of simplicity, we assume henceforth
%that $\Af$ is an elliptic curve, and denote it by~$E$.
%(we expect that most of our discussion should hold with
%$E$ replaced by~$\Af$ below as well).

%Suppose $D \neq -3, -4$ is a negative fundamental discriminant,
%and let $K = \Q(\sqrt{D})$ be the associated quadratic imaginary field.
Let $K$ be a quadratic imaginary field of discriminant not equal
to~$-3$ or~$-4$, and  such that
all primes dividing~$N$
split in~$K$ and such that 
the $L$-function of~$E$ over~$K$ vanishes to order one at~$s=1$.
%${\rm ord}_{s=1} L(E/K, s) =1$.
Choose an ideal~$\NN$ of the ring of integers~$\OO_K$
of~$K$ such that $\OO_K/\NN \isom \Z/N \Z$. Then the complex tori
$\C/\OO_K$ and~$\C/\NN^{-1}$ define elliptic curves related by 
a cyclic $N$-isogeny, hence a complex valued point~$x$ of~$X_0(N)$.
This point, called a Heegner point, is defined over the
Hilbert class field~$H$ of~$K$.
%Let $x$ be a Heegner point of discriminant~$D$ on~$X$ (as in~\cite[\S~I.3]{gross-zagier}). 
Let $P \in J(K)$ be the 
class of the divisor 
$\sum_{\sigma \in {\rm Gal}(H/K)} ((x) - (\infty))^\sigma$.
\comment{
, where
$H$ is the Hilbert class field of~$K$. 

Suppose $D \neq -3, -4$ is a negative fundamental discriminant,
and let $K = \Q(\sqrt{D})$. Let $x$ be a Heegner point of
discriminant~$D$ on~$X$ (as in~\cite[\S~I.3]{gross-zagier}). 
Let $P \in J(K)$ be the 
class of the divisor 
$\sum_{\sigma \in {\rm Gal}(H/K)} ((x) - (\infty))^\sigma$, where
$H$ is the Hilbert class field of~$K$. 
}

By~\cite[p.311--313]{gross-zagier}, the index in~$E(K)$
of the subgroup generated by~$\pi(P)$  is finite; note that
this subgroup is just $\pi(\T P)$.
Also, the order of~$\Sha(\Af/K)$ is finite, by work of Kolyvagin
(see, e.g., \cite[Thm~1.3]{gross:kolyvagin}).
By~\cite[{\S}V.2:(2.2)]{gross-zagier}, the second part of
the Birch and Swinnerton-Dyer conjecture becomes:
\begin{conj}[Birch and Swinnerton-Dyer, Gross-Zagier] \label{conj:bsd}
\begin{eqnarray} \label{gzformula}
|\Af(K)/ \pi(\T P)| \stackrel{?}{=} 
\ce \cdot \prod_{p | N} c_p(\Af) \cdot 
\sqrt{| \Sha(\Af/K) |} \ ,
\end{eqnarray}
where $\ce$ is the Manin constant of~$\Af$ 
(conjectured to be one),
and $c_p(\Af)$ is the Tamagawa number of~$\Af$ at the prime~$p$ 
(i.e., the  order of
the arithmetic component group of~$\Af$ at the prime~$p$).
\end{conj}

%Note that $\pi(\T P) = \Z \pi(P)$ and 
%Note that 
%the index on the left side of~(\ref{gzformula}) 
%is finite %since the analytic rank of~$f$ is one (
The theory of Euler systems can be used to show that the actual
value of the order of~$\Sha(\Af/K)$ divides the order predicted
by the conjectural formula~(\ref{gzformula})
(equivalently, that the right
side of~(\ref{gzformula}) divides the left side), under certain hypotheses,
and staying away from certain primes (see \cite[Theorem~A]{kolyvagin:euler},
and also~\cite[Thm~1.3]{gross:kolyvagin}). 
Our goal is to try to prove results towards
divisibility in the opposite direction, i.e., that the left
side of~(\ref{gzformula}) divides the right side. 
In particular, we shall extract an integer  factor of the left
side of~(\ref{gzformula}) which we will relate to congruences of~$f$ with
eigenforms of analytic rank bigger than one, and these congruences
will in turn be related
to the right side of~(\ref{gzformula}) using the theory of visibility,
under certain hypotheses.

\comment{
suppose $p$ is a prime that divides this index. 
The theory of Euler systems should give results in the opposite
direction for the conjectural equality~(\ref{gzformula}).
\edit{mention them!}
Our goal is to work in the other direction
and try to show that $p$ divides the right hand side.
\edit{we do not know of any other result in the other direction}

If $g$ is a newform of level~$N_g$ dividing~$N$, then
let $S'_g$ denote the subspace of~$S_2(\Gamma_0(N),\C)$
spanned by the forms $g(dz)$ where $d$ ranges over the
divisors of~$N/N_g$. 
Let $I$ be the intersection of the annihilators of~$S'_g$
where $g$ ranges over a subset~$S$ of newforms of level dividing~$N$ that 
have analytic rank at least two and such that $S$ is closed under
complex conjugation.
Let $J'$ denote the quotient 
abelian variety $J/(I_f \cap I)J$. The quotient map $\pi: J \ra J/I_f J$
factors through~$J'= J/(I_f \cap I_g)J$. 
Let $\pi''$ denote the map $J \ra J'$  and $\pi'$ the map $J' \ra \Af$ 
in this factorization. 
Let $B'$ denote the kernel of~$\pi'$.
If $g$ is a newform of level~$N_g$ dividing~$N$, then
let $B_g$ denote the abelian subvariety of~$J_0(N_g)$
associated to~$g$ by Shimura~\cite[Thm.~7.14]{shimura:intro}, 
and let $J_g$ be the sum of
the images of~$B_g$ in~$J=J_0(N)$
%\edit{replace by J everywhere?}
under the usual degeneracy maps.
Then $B'$ is isogenous
to the product of all~$J_g$ where $g$ runs over the elements of~$S$. 
}

Let $T$ be a non-empty set of Galois conjugacy classes of newforms of
level dividing~$N$ and having analytic rank more than one. 
Let $S_T$ denote the
subspace of~$S_2(\Gamma_0(N),\C)$ 
spanned by the forms $g(dz)$, where
$g$ runs over elements in the Galois conjugacy classes in~$T$,
and $d$ ranges over the
divisors of~$N/N_g$, where $N_g$ denote the ``true level'' of~$g$. 
Let $I_T$ denote the annihilator of~$S_T$ under the action of~$\T$.
%Let $I_S = \cap_{g \in S} I_g$. %be the intersection of the annihilators of~$g$
%where $g$ ranges over a subset~$S$ of newforms of level dividing~$N$ that 
%have analytic rank at least two.
Let $J'$ denote the quotient 
abelian variety $J/(I_f \cap I_T)J$. 
The quotient map $\pi: J \ra J/I_f J$
factors through~$J'= J/(I_f \cap I_T)J$. 
Let $\pi'$ denote the map $J' \ra \Af$ and $\pi''$ the map $J \ra J'$ in
this factorization. Let $B'$ denote the kernel of~$\pi'$.
Thus we have the following diagram:
$$\xymatrix{
 &  & J \ar[d]_{\pi''} \ar[rd]^{\pi} & & \\
0 \ar[r] & B' \ar[r] & J' \ar[r]^{\pi'}  & \Af \ar[r] & 0\\
}$$
Note that $J'$ and~$B'$ depend on the choice of the set~$T$; we have
suppressed the dependency in the notation for simplicity
(for certain interesting choices of~$T$, see Section~\ref{sec:rmks}).
Let $\Afp$ denote the image of~$\Afdual \subseteq J$ in~$J'$ under
the  quotient map~\mbox{$\pi'': J \ra J'$} and let $\pi''(\T P)_f$ 
denote the free part of~$\pi''(\T P)$. 
\begin{lem} We have $\pi''(\T P)_f \subseteq \Afp(K)$ with finite index, and
\begin{eqnarray} \label{eqn:fact2}
& & |\Af(K)/ \pi(\T P)|   \\
&  = & \bigg|\frac{J'(K)}{B'(K) + \Afp(K)}\bigg| \cdot 
\big|{\rm ker}\big(H^1(K,B') \ra H^1(K,J')\big) \big|
\cdot  \frac{\big|\frac{B'(K) + \Afp(K)}{B'(K) + \pi''(\T P)_f}\big|}
{\big|\frac{B'(K) + \pi''(\T P)}{B'(K) + \pi''(\T P)_f}\big|} \ \ \nonumber.
\end{eqnarray}
\end{lem}
\begin{proof}
By~\cite[Prop.~1.6]{ag:visrk1}, we have 
\begin{eqnarray} \label{eqn:fact}
|\Af(K)/ \pi(\T P)| 
 =  \bigg|\frac{J'(K)}{B'(K) + \pi''(\T P)}\bigg| \cdot 
\big|{\rm ker}\big(H^1(K,B') \ra H^1(K,J')\big)\big|.
\end{eqnarray}
%Suppose henceforth that $r$ is coprime to $E(K)_{\rm tor}$. Then 
If $h$ is an eigenform of weight~$2$ on~$\Gamma_0(N)$,
then $\T P \cap \Ahdual(K)$ is infinite if and only if $h$ has
analytic rank one (this follows by~\cite[Thm~6.3]{gross-zagier}
if $h$ has analytic rank bigger than one,
and the fact that $\Ahdual(K)$
is finite if~$h$ has analytic rank zero, by~\cite{kollog:finiteness}).
The composite $\Afdual \stackrel{\pi''}{\ra} \Afp \stackrel{\pi'}{\ra} \Af$
is an isogeny, and so
$J'$ is isogenous to $\Afp \oplus B'$. 
Now $B'$ is isogenous
to a product of~$\Agdual$'s (with multiplicities)
where $g$ runs over representatives
of conjugacy classes of eigenforms in~$T$; all these 
eigenforms have analytic rank
greater than one.  
Thus from the discussion above,  we 
see that the free part of $\Afp(K)$ contains $\pi''(\T P)_f$. 
%Let us denote the free part of~$\pi''(\T P)$ by~$\pi''(\T P)_f$; so
%$\pi''(\T P)_f \subseteq \Afp(K)$.
%We can now rewrite equation~(\ref{eqn:fact}) as
The lemma now follows from equation~(\ref{eqn:fact}).
We remark that the transition from equation~(\ref{eqn:fact}) to
equation~(\ref{eqn:fact2}) is
analogous to the situation
in the rank one case (cf. Theorem~3.1 of~\cite{agmer}
and its proof), where the idea is due to L.~Merel.
\end{proof}

\comment{
The reason for factoring out the product
$$\bigg|\frac{J'(K)}{B'(K) + \Afp(K)}\bigg|\cdot 
\big|{\rm ker}\big(H^1(K,B') \ra H^1(K,J')\big) \big|$$
in equation~(\ref{eqn:fact2}) from the quantity
$|\Af(K)/ \pi(\T P)|$, which is the left side of
the Birch and Swinnerton-Dyer conjectural formula~(\ref{gzformula}),
is that 
we can say something about it:
}

The reason for factoring the quantity
$|\Af(K)/ \pi(\T P)|$, which is the left side of
the Birch and Swinnerton-Dyer conjectural formula~(\ref{gzformula}),
as above in equation~(\ref{eqn:fact2}) is that we can 
say something about the factor
$$\bigg|\frac{J'(K)}{B'(K) + \Afp(K)}\bigg|\cdot 
\big|{\rm ker}\big(H^1(K,B') \ra H^1(K,J')\big) \big|$$
in this factorization:

\begin{prop} \label{prop:main}
Suppose $q$ is a prime that 
%does not divide the order of the torsion part of~$\pi''(\T P)$, and 
divides the product
$$\bigg|\frac{J'(K)}{B'(K) + \Afp(K)}\bigg|\cdot 
\big|{\rm ker}\big(H^1(K,B') \ra H^1(K,J')\big) \big|.$$ 
Then $q$ divides the order of~$B' \cap \Afp$, and 
there is an eigenform~$g$ in~$S_T$
%the subspace of~$S_2(\Gamma_0(N),\C)$
%spanned by the representative newforms in~$T$
%on~$\Gamma_0(N)$ having
%analytic rank greater than one
such that $f$ is congruent to~$g$ modulo a prime ideal $\qq$
over~$q$ in the ring of integers 
of the number field generated by the Fourier coefficients
of~$f$ and~$g$. 
\end{prop}
We will prove this Proposition in Section~\ref{sec:proofs}.
If $G$ is a finite abelian group and $r$ is a prime, then let $|G|_r$ denote
the order of the $r$-primary component of~$G$ (equivalently,
$|G|_r$ is the highest power of~$r$ that divides the order of~$G$).
If $q$ is a prime that does not divide the order of the torsion subgroup
of~$\pi''(\T P)$, then 
in view of equation~(\ref{eqn:fact2}),
the Birch and Swinnerton-Dyer conjectural formula~(\ref{gzformula})
says:
\begin{eqnarray} \label{eqn:fact3}
& & \big|\caf \cdot \prod_{p | N} c_p(\Af)\big|_q \cdot 
\sqrt{| \Sha(\Af/K) |_q}  \\
& \stackrel{?}{=} & \bigg|\frac{J'(K)}{B'(K) + \Afp(K)}\bigg|_q \cdot 
\big|{\rm ker}\big(H^1(K,B') \ra H^1(K,J')\big) \big|_q
\cdot  \bigg|\frac{B'(K) + \Afp(K)}{B'(K) + \pi''(\T P)_f}\bigg|_q. \nonumber
\end{eqnarray}
In particular, the conjecture predicts that
the product 
$$\bigg|\frac{J'(K)}{B'(K) + \Afp(K)}\bigg|_q \cdot 
\big|{\rm ker}\big(H^1(K,B') \ra H^1(K,J')\big) \big|_q$$ divides 
$c_{\scriptscriptstyle E} \cdot \prod_{p | N} c_p(E) \cdot 
\sqrt{| \Sha(E/K)|}$.
\comment{
Since we are assuming that 
$r$ is coprime to the order of the torsion part of~$\pi''(\T P)$
(which is condition~(c) on~$r$),
considering that $r/r'$ divides
$|\frac{J'(K)}{B'(K) + \Afp(K)}|$,
we see that $r/r'$ divides
$|\frac{J'(K)}{B'(K) + \pi''(\T P)}|$ as well.

Let $B = {\rm ker\ \pi}$, which is an abelian subvariety of~$J_0(N)$.
If $g$ is an eigenform of weight~$2$ on~$\Gamma_0(N)$, then
we shall denote the dual abelian variety of~$A_g$ by~$\Agdual$; it
is an abelian subvariety of~$J_0(N)$. 
If $H$ is a subgroup of a finitely-generated abelian group~$G$,
then the {\em saturation} $H$ in~$G$ is the largest subgroup 
of~$G$ containing~$H$ with finite index. Let $\widehat{\T P}$ denote
the saturation of~$\T P$ in~$J(K)$. The following result expresses
the integer~$|E(K)/ \pi(\T P)|$ on the left side of~(\ref{gzformula})
as a product of three integers:

\begin{prop} \label{prop:factorization}
We have
$$|E(K)/ \pi(\T P)| = 
\bigg|\frac{J(K)}{B(K) + \tphat}\bigg| \cdot 
\bigg|\frac{B(K) + \tphat}{B(K) + \T P}\bigg| \cdot
|{\rm ker}\big(H^1(K,B) \ra H^1(K,J)\big)|.$$
\end{prop}

We give the proof of this proposition in Section~\ref{sec:proofofprop}.
Let $D$ denote the abelian subvariety of~$J_0(N)$ generated by~$\Agdual$
where $g$ ranges over all eigenforms other than~$f$
of weight~$2$ on~$\Gamma_0(N)$ and
having analytic rank one.
}

Using Proposition~\ref{prop:main} and the theory of visibility, we 
can show the following result towards this predicted divisibility:

\begin{thm} \label{thm:main}
Let  $q$ be a prime that 
%does not divide the order of the torsion part of~$\pi''(\T P)$, and 
divides the product
$$\bigg|\frac{J'(K)}{B'(K) + \Afp(K)}\bigg|\cdot 
\big|{\rm ker}\big(H^1(K,B') \ra H^1(K,J')\big) \big|.$$ 
Suppose that
$q \ndiv 2 N$ and that $E[q]$
%for all maximal ideals~$\qq$ of~$\T$
%with residue characteristic~$q$, $\Af[\qq]$ 
is an irreducible
representation of the absolute Galois group of~$K$.
%and $A_g[\qq]$ is irreducible.
Assume that the first part of the Birch and Swinnerton-Dyer conjecture
holds for all quotients of~$J_0(N)$ associated to eigenforms
of analytic rank greater than one.
Suppose that  
for all primes $p \divs N$, $p \not\equiv - w_p \pmod q$,
where $w_p$ is the sign of the Atkin-Lehner involution~$W_p$ acting on~$f$.
We have two possibilities:

%Consider the following statement:\\
%Suppose~(*) is false. 

\noindent
Case (i) For all primes $p \divs N$,
$f$ is not congruent to a newform~$g$
of level dividing $N/p$ 
(for Fourier coefficients of index coprime to~$Nq$)
modulo a prime ideal over~$q$ in the ring of integers of the number
generated by the coefficients of~$f$ and~$g$:\\
In this case, suppose that for all primes~$p$ such that $p^2 \divs N$, 
we have $p \not\equiv -1  \pmod q$.  \\
Then $q$ divides~$\Mid \Sha(\Af) \miD$.

%If for some prime $p \divs N$,
%$f$ is congruent $\bmod \ \qq$ to a newform
%of level dividing $N/p$ (for Fourier coefficients of index coprime to~$Nq$),
%then $\qq$ divides $c_p(\Af)$ provided the following statement is false: \\

\noindent Case (ii) For some prime $p$ dividing~$N$, 
$f$ is congruent to a newform~$g$
of level dividing $N/p$ 
(for Fourier coefficients of index coprime to~$Nq$)
modulo a prime ideal~$\qq$ over~$q$ in the ring of integers of the number
generated by the coefficients of~$f$ and~$g$:\\
In this case, suppose 
that there a triple $p$, $g$, and~$\qq$ as above 
%a prime~$p$ dividing~$N$
%and a newform~$g$ of level dividing $N/p$ such that $f$ is congruent to~$g$
%(for Fourier coefficients of index coprime to~$Nq$)
%modulo a prime ideal over~$q$ in the ring of integers of the number
%generated by the coefficients of~$f$ and~$g$, 
such that $p^2 \ndiv N$, $w_p=-1$, and $A_g[\qq]$ irreducible
as a representation of the absolute Galois group of~$\Q$.\\
%Consider the following statement:\\
%\noindent (*) for all primes $p \divs N$, 
%%$A[\qq]$ is unramified, then $w_p=1$.\\
%if $f$ is congruent $\bmod \ \qq$ 
%to a newform~$h$ of level dividing $N/p$ 
%(for Fourier coefficients of index coprime to~$Nq$), then either
%$p^2 \divs N$, or $w_p=1$, or
%$A_h[\qq]$ is reducible.\\
%If (*) is false, 
Then $q$ divides $\prod_{\scriptscriptstyle{{p\mid N}}} c_p(\Af)$.

\comment{
Suppose that $q \ndiv 2 N$, and for all maximal ideals~$\qq$ of~$\T$
with residue characteristic~$q$, $\Af[\qq]$ is irreducible.
Assume the first part of the Birch and Swinnerton-Dyer conjecture
for all newforms of level dividing~$N$ and having analytic rank
at least one.
Suppose that  
for all primes $p \divs N$, $p \not\equiv - w_p \pmod q$ 
and $p \not\equiv -1  \pmod q$ if $p^2 \divs N$.  
We have two possibilities:\\
\noindent
Case (i) For all primes $p \divs N$,
$f$ is not congruent modulo~$\qq$ to a newform
of level dividing $N/p$ (for Fourier coefficients of index coprime to~$Nq$):\\
In this case, $q^2$ divides $\Mid \Sha(\Af) \miD$. 

\noindent Case (ii) $f$ is congruent modulo~$\qq$ to a newform of lower level
(for Fourier coefficients of index coprime to~$Nq$):\\
In this case, suppose that \\
(*) there is a prime $p$ dividing~$N$ 
such that $p^2 \ndiv N$, $w_p=-1$, 
$A_h[\qq]$ is irreducible, and 
$f$ is congruent modulo~$\qq$ 
to a newform~$h$ of level dividing $N/p$ 
(for Fourier coefficients of index coprime to~$Nq$). \\
Then $q$ divides $\prod_{\scriptscriptstyle{{p\mid N}}} c_p(\Af)$.
}
\end{thm}
\begin{proof}
%First suppose that 
%$f$ is not congruent modulo a prime over~$q$ to a newform of 
%lower level, 
%we see that $g$ is a newform. 
By Proposition~\ref{prop:main},
there is an eigenform~$g$ on~$\Gamma_0(N)$ having
analytic rank greater than one
such that $f$ is congruent to~$g$ modulo a prime ideal $\qq$
over~$q$ in the ring of integers~$\OO$
of the number field generated by the Fourier coefficients
of~$f$ and~$g$. 

Case~(i) now follows from  Theorem~6.1 of~\cite{dsw}, as we now indicate
(for details of some of the definitions below, see~\cite{dsw}).
%Since $q$ is odd,
%considering that the eigenvalue of the Atkin-Lehner involution~$W_N$
%on~$g$ is the same as that of~$f$,
%we see that $g$ has odd analytic rank, and hence its analytic 
%rank is at least~$3$.
%Since we are assuming the first part of
%the Birch and Swinnerton-Dyer conjecture, this implies
%that $A_g$ has Mordell-Weil rank at least $3 \cdot \dim(\Af) = 3$.
Let $T_\qq$ denote the $\qq$-adic Tate module of~$\Afdual = E$.  Let $L_\qq$
denote the quotient field of~$\OO_\qq$, let 
$V_\qq = T_\qq \tensor_{\OO_\qq} L_\qq$, and let $A_\qq$ denote $V_\qq/ T_\qq$.
We denote the corresponding objects for $\Agdual$ by~$T'_\qq$, $V'_\qq$,
and~$A'_\qq$. 
Let $r$ denote the dimension of~$H^1_f(K, V'_\qq)$ over~$L_\qq$.
Then $r$ is at least the analytic rank of~$g$ (since we are assuming
the first part of the Birch and Swinnerton-Dyer conjecture for~$\Agdual$),
i.e., at least~$2$. Theorem~6.1 of loc. cit. (which is stated over~$\Q$,
but works over~$K$ as well) tells us that 
%Since the other hypotheses of Theorem~6.1 of loc. cit.
%are satisfied, the conclusion of  then tells us that
the $\qq$-torsion subgroup of the 
Selmer group~$H^1_f(K, A_\qq)$ of~$\Edual$ 
has~$\OO_\qq/\qq$-rank at least~$r$. Since 
%under our assumption, $\Edual(K)[q] = 0$, 
%\edit{At this point you don't need to use the vanishing
%of the rational q-torsion. Also, it only makes sense to
%talk of $F_q$-ranks of q-torsion subgroups.
%(make it the $F_q$ rank of the q-torsion subgroup of ...).}
the abelian group~$\Edual(K)$ has rank one, 
the image of $H^1_f(K, V_\qq)$ in the 
$\qq$-torsion subgroup of~$H^1_f(K, A_\qq)$  
has $\OO_\qq/\qq$-rank at most one. This shows that
$|\Sha(\Edual/K)|$ is divisible by $q^{r-1}$, in particular 
by~$q^{2-1}=q$ (since $r \geq 2$).
By the perfectness of the Cassels-Tate
pairing, we see that $q$ divides the order of~$\Sha(E/K)$ as well.
This proves Case~(i).
Case~(ii) follows from~\cite[Prop.~6.3]{agmer}.
\end{proof}

\begin{cor}
Let $q$ be a prime that 
%does not divide the order of the torsion part of~$\pi''(\T P)$, and 
divides the product
$$\bigg|\frac{J'(K)}{B'(K) + \Afp(K)}\bigg| \cdot 
\big|{\rm ker}\big(H^1(K,B') \ra H^1(K,J')\big) \big|.$$ 
Suppose that $N$ is prime, $q \ndiv N(N+1)$, and $E[q]$ is an
irreducible representation of the absolute Galois group of~$K$.
Assume that the first part of the Birch and Swinnerton-Dyer conjecture 
holds for all quotients of~$J_0(N)$ associated to eigenforms
of analytic rank
greater than one. 
%Assume that for all newforms~$g$ on~$\Gamma_0(N)$,
%if $L(g,1) = 0$, then
%the rank of~$A_g(\Q)$ is positive (which would follow if
%the first part of the Birch and Swinnerton-Dyer conjecture were true).
Then $q$ divides $\Mid \Sha(\Af) \miD$.
\end{cor}
\begin{proof}
%By~\cite[Prop.~14.2]{mazur:eisenstein}, 
%since $q \ndiv 2 \cdot {\rm numr}(\frac{N-1}{12})$, 
%it follows that $\Af^\vee[\qq]$ is  irreducible
%for all maximal ideals~$\qq$ of~$\T$
%with residue characteristic~$q$.
%for any newform~$h$ in~$S_2(\Gamma_0(N), \C)$ (in particular, if $h=f$ or $h=g$). 
Since $w_N = 1$, we have $N \not\equiv - 1 \pmod q$ by hypothesis.
Also, since the level is prime, there are no newforms of lower level.
The corollary now follows from Theorem~\ref{thm:main}.
\end{proof}

We remark that N.~Dummigan has informed us that the hypothesis
that for all 
primes~$p$ dividing~$N$,
$p \not\equiv -w_p \bmod q$ can be eliminated from~\cite[Thm.~6.1]{dsw},
and hence from Theorem~\ref{thm:main}; if this is the case, then
the hypothesis that $q \ndiv (N+1)$ can be eliminated from the corollary
above.

\comment{As mentioned earlier, 
if $q$ is a prime that does not divide the order of the torsion part
of~$\pi''(\T P)$, 
the second part of the Birch and Swinnerton-Dyer
conjecture~(\ref{gzformula}) predicts that
the product 
$$\bigg|\frac{J'(K)}{B'(K) + \Afp(K)}\bigg|_q \cdot 
\big|{\rm ker}\big(H^1(K,B') \ra H^1(K,J')\big) \big|_q$$ divides 
$c_{\scriptscriptstyle E} \cdot \prod_{p | N} c_p(E) \cdot 
\Mid \Sha(E/K) \miD^{1/2}$.
}

In view of our discussion just preceding Theorem~\ref{thm:main},
the theorem and corollary above are partial results towards
the second part of
the Birch and Swinnerton-Dyer conjecture in the analytic rank one case,
and provide theoretical evidence supporting the conjecture.
Also, under certain hypotheses (the most serious of which is
the first part of the Birch and Swinnerton-Dyer conjecture),
we have shown that if a prime~$q$ divides a certain factor of the left
side of the Birch and Swinnerton-Dyer conjectural formula~(\ref{gzformula}), then $q$ divides
the right side the formula (which includes~$\sqrt{|\Sha(\Af)|}$ as a factor).
Thus our result is a step in trying to prove that 
the left side of the Birch and Swinnerton-Dyer conjectural formula~(\ref{gzformula})
divides the right side.
As mentioned earlier, 
the theory of Euler systems gives results in the opposite direction,
viz., that the right side of the Birch and Swinnerton-Dyer conjectural formula~(\ref{gzformula})
divides the left side (under certain hypotheses).
Thus our result
fits well in the ultimate goal of trying to prove the second part
of the Birch and Swinnerton-Dyer conjecture in the analytic rank one case. Note
that the theory of Euler systems can also be used to 
construct non-trivial elements of the Shafarevich-Tate group
(e.g., see~\cite{mccallum:kolyvagin}).

In Section~\ref{sec:rmks} we make some further remarks about our main result
and in Section~\ref{sec:proofs}, we 
give the proof of Proposition~\ref{prop:main}. \\

\noindent {\it Acknowledgements:} We are grateful to Neil Dummigan for answering
some questions regarding~\cite{dsw}. We would also like to thank John Cremona
and Mark Watkins for some numerical data that encouraged the author to pursue
the investigations in this article.

\section{Some further remarks} \label{sec:rmks}

Note that the product 
$$\bigg|\frac{J'(K)}{B'(K) + \Afp(K)}\bigg|\cdot 
\big|{\rm ker}\big(H^1(K,B') \ra H^1(K,J')\big) \big|$$
depends on the choice of~$T$ in Section~\ref{sec:intro}.
There are two interesting choices of~$T$, which are the two extreme cases.

The first is where $T$ consists of conjugacy classes of {\em all} newforms
of level dividing~$N$ that have analytic rank more than one. Then
by Proposition~\ref{prop:main}, the product above is allowed to be
divisible by all primes~$p$ such that 
$f$ is congruent to an eigenform~$g$ of analytic rank bigger
than one modulo some prime
ideal over~$p$; such eigenforms~$g$ are precisely all of
the eigenforms on~$\Gamma_0(N)$
using which the theory
of visibility can be used to construct non-trivial elements of the
the Shafarevich-Tate group of~$E$ as in Theorem~\ref{thm:main}. Thus in some sense,
we are getting the most out of the theorem with this choice of~$T$.

The other extreme choice of~$T$ is that it consists of the conjugacy class
of a {\em single} newform on~$\Gamma_0(N)$
that has analytic rank more than one. The advantage of this choice is
that we are able to prove a sort of converse to Proposition~\ref{prop:main}:

%The following result is easily extracted from our companion 
%paper~\cite{ag:visrk1}.
\begin{prop} \label{prop:revdiv}
Recall that $f$ is a newform with
integer Fourier coefficients that has analytic rank one. Suppose  
there is a newform~$g$ with integral Fourier coefficients that has 
analytic rank greater than one
such that 
$f$ and~$g$ are congruent modulo an odd prime~$q$.
Take $T$ to be
the singleton set~$\{g\}$ in the definition of $J'$ and~$B'$
in Section~\ref{sec:intro}. Suppose that
$q^2 \nmid N$, $A^{\vee}_f[q]$ and~$\Agdual[q]$ are irreducible
representations of the absolute Galois group of~$\Q$, and
$q$ does not divide the order of the torsion subgroup 
of the projection of~$\T P$ in~$J'$.
Then $q$ divides the product
$$\bigg|\frac{J'(K)}{B'(K) + \Afp(K)}\bigg|\cdot 
\big|{\rm ker}\big(H^1(K,B') \ra H^1(K,J')\big) \big|.$$
\end{prop}
\begin{proof}
%Then by Corollary~1.5 and Proposition~1.6 of~\cite[1.6]{ag:visrk1} 
The proposition essentially follows from the main result of~\cite{ag:visrk1},
as we now indicate.
Take $E=A^{\vee}_f$ and $F=\Agdual$ in loc. cit.
Then since $B'$ is isogenous to ~$\Agdual$, we see that $F' = B'$
in the notation of loc. cit.
Let $r$ denote the highest power of~$q$ modulo which $f$ and~$g$
are congruent. Then condition (a) on~$r$
in Theorem~1.4 of loc. cit. is satisfied.
By the discussion just before Lemma 1.2 of loc. cit.,
the hypotheses that $q$ is odd, $q^2 \nmid N$, and 
$A^\vee_f[q]$ and~$\Agdual[q]$ are  irreducible 
imply that $r$ satisfies condition~(b) in Theorem~1.4 of loc. cit.
The hypothesis that 
$q$ does not divide the order of the torsion subgroup 
of the projection of~$\T P$ in~$J'$ 
implies that condition~(c) 
in Theorem~1.4 of loc. cit. is satisfied. 
Then the proof of Theorem~1.4 of loc. cit. shows that
$r$ divides
$$\bigg|\frac{J'(K)}{B'(K) + \Afp(K)}\bigg|\cdot 
\big|{\rm ker}\big(H^1(K,B') \ra H^1(K,J')\big) \big|$$
(the statement just after Lemma 2.1 of loc. cit. shows that 
$r/r'$ divides the first factor in the product
above, and the very last statement
shows that $r'$ divides the second factor in the product above, 
where $r'$ is a 
certain divisor of~$r$ defined in loc. cit.). The proposition now follows.
\end{proof}

%Then as mentioned in the proof of 
As mentioned before, the proposition above is a result that is in a direction 
opposite  to that of Proposition~\ref{prop:main},
and is a partial result in trying to characterize the primes
that divide the factor
$$\bigg|\frac{J'(K)}{B'(K) + \Afp(K)}\bigg|\cdot 
\big|{\rm ker}\big(H^1(K,B') \ra H^1(K,J')\big) \big|$$
of the ``analytic'' left side of 
the Birch and Swinnerton-Dyer conjectural formula~(\ref{gzformula}), which we
related to the ``arithmetic'' right side of this formula.
Notice the similarity with the rank zero case in~\cite{agmer},
where we isolated a factor of 
the ``analytic'' left side of the Birch and Swinnerton-Dyer formula
that could 
be characterized in terms of congruences analogous to the ones above
and related these congruences to the ``arithmetic'' right side (the results for
the analytic rank zero case are more precise).

\begin{rmk} \label{rmk:one}
One question that remains is whether
the product 
$$\bigg|\frac{J'(K)}{B'(K) + \Afp(K)}\bigg|\cdot 
\big|{\rm ker}\big(H^1(K,B') \ra H^1(K,J')\big) \big|$$ 
is non-trivial in general, and if so, how often.
It would be nice to have some numerical data where the hypotheses
of Proposition~\ref{prop:revdiv} are satisfied, 
%when the situation described two  paragraphs above does take
%place, 
so that the product above
is non-trivial. If this happens, then in view of
Theorem~\ref{thm:main}, we expect that either~$\Sha(\Af)$
is non-trivial or an arithmetic component group of~$\Af$ is non-trivial,
of which the former seems more likely.
Since it is difficult to compute the actual order of the Shafarevich-Tate group, 
we looked at the Birch and Swinnerton-Dyer conjectural
orders in Cremona's online ``Elliptic curve data''~\cite{cremona:online}.
Unfortunately the conjectural orders of the Shafarevich-Tate groups of elliptic
curves of analytic rank one at low levels are usually one or powers of~$2$,
which makes it difficult to find
examples where the hypotheses of Proposition~\ref{prop:revdiv} 
can be verified easily. 
For levels up to~$30000$, we found only
one optimal elliptic curve of Mordell-Weil rank one for which
the conjectural order of the Shafarevich-Tate group was divisible
by an {\em odd} prime: the curve~$E$
with label 28042A, for which the
conjectural order of the Shafarevich-Tate group is~$9$. 
At the same level, the curve $F=\ $28042B has Mordell-Weil rank~$3$ and 
the newforms~$f$ and~$g$  corresponding to 28042A and 28042B 
respectively have Fourier
coefficients that are congruent modulo~$3$ for every prime index up to~$100$
(although this is not enough to conclude that the newforms are congruent
modulo~$3$ for all Fourier coefficients, 
cf.~\cite{agashe-stein:schoof-appendix}). 
We do not know how to verify
the hypotheses in Proposition~\ref{prop:revdiv}
that 
$3$ does not divide the order of the torsion subgroup 
of the projection of~$\T P$ in~$J'$ and that 
$E[3]$ and~$F[3]$ are irreducible representations
of the Galois group of~$K$ 
(we remark though
that  by~\cite{cremona:online}, $E$ and~$F$ have no $3$-torsion over~$\Q$).
So while we  cannot be sure that Proposition~\ref{prop:revdiv} applies
to show that 
%are congruent modulo~$3$, 
$3$ divides the product
$$\bigg|\frac{J'(K)}{B'(K) + \Afp(K)}\bigg|\cdot 
\big|{\rm ker}\big(H^1(K,B') \ra H^1(K,J')\big) \big|,$$
it is quite encouraging that for the first example where 
the conjectural order of the Shafarevich-Tate group of an elliptic
curve is divisible by an odd prime, there is a congruence modulo the
same  prime that 
might show that the product above is divisible by the prime in question, 
and hence explain why the prime divides the order 
of the Shafarevich-Tate group. 
\end{rmk}

\begin{rmk} \label{rmk:two}
Let $q$ be a prime that
does not divide the order of the torsion subgroup
of~$\pi''(\T P)$. 
Recall the conjectural equation~(\ref{eqn:fact3}), which is predicted by the
Birch and Swinnerton-Dyer conjecture, and which we repeat below:
\begin{eqnarray} 
& & \big|\caf \cdot \prod_{p | N} c_p(\Af)\big|_q \cdot 
\sqrt{ | \Sha(\Af/K) |_q} \nonumber \\
& \stackrel{?}{=} & \bigg|\frac{J'(K)}{B'(K) + \Afp(K)}\bigg|_q \cdot 
\big|{\rm ker}\big(H^1(K,B') \ra H^1(K,J')\big) \big|_q
\cdot  \bigg|\frac{B'(K) + \Afp(K)}{B'(K) + \pi''(\T P)}\bigg|_q. \nonumber  
\end{eqnarray}
While we were able to relate certain primes dividing the product
$$\bigg|\frac{J'(K)}{B'(K) + \Afp(K)}\bigg|_q\cdot 
\big|{\rm ker}\big(H^1(K,B') \ra H^1(K,J')\big) \big|_q$$
on the right side of the equation above to 
its left side, 
one question that remains is to interpret
the remaining
factor $$\bigg|\frac{B'(K) + \Afp(K)}{B'(K) + \pi''(\T P)}\bigg|_q$$
on the right side of equation, which is expected to divide
the left side of equation above.
%, which is $\big|\caf \cdot \prod_{p | N} c_p(\Af)\big|_q \cdot 
%\Mid \Sha(\Af/K) \miD^{1/2}_q$.
Now we were able to relate the primes dividing the product
$$\bigg|\frac{J'(K)}{B'(K) + \Afp(K)}\bigg|\cdot 
\big|{\rm ker}\big(H^1(K,B') \ra H^1(K,J')\big) \big|$$
to $\prod_{p | N} c_p(\Af) \cdot \sqrt{ | \Sha(\Af/K) |}$ in
Theorem~\ref{thm:main} using the theory of visibility
and the existence of congruences modulo prime ideals over~$q$
with eigenforms at the {\em same} level~$N$
that have analytic rank more than one. If $M$ is a positive integer,
then $f$ can be mapped to~$S_2(\Gamma_0(NM), \C)$ using suitable
degeneracy maps, and if there is an eigenform at the {\em higher} level~$NM$ that
is congruent to the image of~$f$ in~$S_2(\Gamma_0(NM), \C)$
modulo some prime ideal over a prime~$q$, then again the theory
of visibility can sometimes be used to show that $q$ divides
the order of~$\Sha(E)$ (e.g., see~\cite[\S4.2]{agst:vis}). We loosely call
this phenomenon ``visibility at higher level''.
It has been conjectured that any element 
the Shafarevich-Tate group can be explained by 
visibility at some higher level 
(see~Conjecture~7.1.1 in~\cite{jetchev-stein}
for details and a precise statement).
Thus we suspect that one may be able to explain the factor 
$$\bigg|\frac{B'(K) + \Afp(K)}{B'(K) + \pi''(\T P)}\bigg|_q$$
using the idea of visibility at higher level, at least in specific examples.
The situation is 
%visibility via  congruences of~$f$ with newforms at level a multiple of~$N$
%and having analytic rank bigger than one (cf.~\cite[Remark~6.10(3)]{agmer}).
%The rationale is that this
similar to the case where $f$ has analytic rank one~\cite{agmer}, when
we were able to understand a certain factor using the theory
of visibility and congruences of~$f$ with eigenforms of higher rank on~$\Gamma_0(N)$,
and suspected that to explain the remaining factor, one would need
to use visibility at a higher level.
\end{rmk}

In view of Remarks~\ref{rmk:one} and~\ref{rmk:two},
we hope that our article motivates more detailed
computations similar to those in~\cite{agst:bsd}
for the analytic rank one case,
especially since all this pertains to the Birch and Swinnerton-Dyer
conjecture.

\section{Proof of Proposition~\ref{prop:main}}
\label{sec:proofs}

%We start by giving the proof of Theorem~\ref{thm:main}, which
%depends only on Proposition~\ref{prop:main}. After that, we give
%the proof of Proposition~\ref{prop:main}.
%The rest of this section is devoted to the proof of
%Proposition~\ref{prop:main}.
Following a similar situation in~\cite{cremona-mazur}, 
consider the short exact sequence
\begin{eqnarray} \label{eqn:ses}
0 \ra B' \cap \Afp \ra B' \oplus \Afp \ra J' \ra 0,
\end{eqnarray}
where the map $B' \cap \Afp \ra B' \oplus \Afp$
is the anti-diagonal embedding $x \mapsto (-x,x)$ and the map 
$B' \oplus \Afp \ra J'$ is given by $(b',e') \mapsto b'+e'$.

\begin{lem} \label{lem1}
Suppose $q$ is a prime that 
divides 
$$\bigg|\frac{J'(K)}{B'(K) + \Afp(K)}\bigg|.$$
Then $q$ divides the order of~$B' \cap \Afp$.
\end{lem}
\begin{proof}
The long exact sequence associated to~(\ref{eqn:ses}) gives us
$$\cdots \ra B'(K) \oplus \Afp(K) \ra J'(K) \ra 
H^1(K,B' \cap \Afp) \ra H^1(K,B' \oplus \Afp) \ra \cdots,$$
from which we get
\begin{eqnarray} \label{galcoh}
\frac{J'(K)}{B'(K) + \Afp(K)} = 
\ker\big(H^1(K,B' \cap \Afp) \ra H^1(K,B' \oplus \Afp)\big).
\end{eqnarray}

Since $q$ divides
$|\frac{J'(K)}{B'(K) + \Afp(K)}|$, there is 
an element~$\sigma$ of the right hand side of~(\ref{galcoh})
of order~$q$. Since $B' \cap \Afp$ is finite, so is $H^1(K,B' \cap \Afp)$,
and its order divides $|B' \cap \Afp|$. Hence $q$ divides $|B' \cap \Afp|$. 
\end{proof}

\begin{lem}  \label{lem2}
Suppose $q$ is a prime that 
divides $\big|{\rm ker}\big(H^1(K,B') \ra H^1(K,J')\big) \big|$.
Then $q$ divides the order of~$B' \cap \Afp$.
\end{lem}
\begin{proof}
By hypothesis, there is an element~$\sigma$
of~${\rm ker}\big(H^1(K,B') \ra H^1(K,J')\big)$ of order~$q$.
The long exact sequence associated to~(\ref{eqn:ses}) gives us
\begin{eqnarray} \label{eqn:les2}
\cdots \ra H^1(K, B' \cap \Afp) \ra H^1(K, B') \oplus H^1(K,\Afp) 
\ra H^1(K,J') \ra \cdots.
\end{eqnarray}
The element $(\sigma, 0) \in H^1(K,B') \oplus H^1(K,\Afp)$ of order~$q$
in the middle group in~(\ref{eqn:les2})
maps to zero in the rightmost group~$H^1(K,J')$ in~(\ref{eqn:les2}),  
and thus by the exactness of~(\ref{eqn:les2}),
there is a non-trivial element $\sigma' \in H^1(K,B' \cap \Afp)$ of order
divisible by~$q$
that maps to $(0, \sigma) \in H^1(K,B') \oplus H^1(K,\Afp)$.
Again, since $B' \cap \Afp$ is finite, so is $H^1(K,B' \cap \Afp)$,
and its order divides $|B' \cap \Afp|$. Hence $q$ divides $|B' \cap \Afp|$. 
\end{proof}

\begin{proof}[Proof of Proposition~\ref{prop:main}]
By Lemmas~\ref{lem1} and~\ref{lem2}, 
we see that $q$ divides the order of~$B' \cap \Afp$, which proves the first
claim in Proposition~\ref{prop:main}. 
The second claim follows from the first, using an argument similar to
the one in~\cite[\S5]{agmer}, which in turn mimics the proof
that the modular degree divides the congruence number~\cite{ars:moddeg},
as we now explain.
%The second claim
%follows from the first by a generalization of the proof 
%of Theorem~3.6(a) in~\cite{ars:moddeg}, which says that 
%the modular exponent divides the congruence exponent,
%as we now indicate.

If $h$ is a newform of level~$N_h$ dividing~$N$, then
let $B_h$ denote the abelian subvariety of~$J_0(N_h)$
associated to~$h$ by Shimura~\cite[Thm.~7.14]{shimura:intro}, 
and let $J_h$ denote the sum of
the images of~$B_h$ in~$J=J_0(N)$
%\edit{replace by J everywhere?}
under the usual degeneracy maps; 
note that $J_h$ depends only on the Galois conjugacy class of~$h$.
Let $C$ denote $(I_f \cap I_T)J$.
Then $C$ is the abelian subvariety of~$J$ generated
by~$J_g$ where $g$ ranges over Galois conjugacy classes of newforms
of level dividing~$N$
other than orbit of~$f$ and other than the
classes in~$T$.
%and such that $g$ has analytic rank at most one; 
Let $B$ denote abelian subvariety of~$J$ generated
by~$J_g$ where $g$ ranges over Galois conjugacy classes of newforms
of level dividing~$N$
other than the orbit of~$f$
and let $A$ denote abelian subvariety of~$J$ generated
by~$J_g$ where $g$ ranges over Galois conjugacy classes of newforms
of level dividing~$N$
other than the
classes in~$T$.
Then $\Afp = A/C$ and $B'= B/C$.
Now applying the arguments of~\cite[\S5]{agmer} but with $A$, $B$, and~$C$
as above, 
the fact that $q$ divides
the order of~$E' \cap B' = A/C \cap B/C$ implies that
there is an eigenform~$g$ in the subspace of~$S_2(\Gamma_0(N),\C)$
generated by the newforms in~$T$ such that 
$f$ is congruent to~$g$ modulo a prime ideal $\qq$
over~$q$  in the ring of integers
of the number field generated by the Fourier coefficients
of~$f$ and~$g$.
The second claim in Proposition~\ref{prop:main} follows, considering
that every newform~$g$ in~$T$ has analytic rank more than one.
\end{proof}

\bibliographystyle{amsalpha}         

\providecommand{\bysame}{\leavevmode\hbox to3em{\hrulefill}\thinspace}
\providecommand{\MR}{\relax\ifhmode\unskip\space\fi MR }
\providecommand{\MRhref}[2]{%
  \href{http://www.ams.org/mathscinet-getitem?mr=#1}{#2}
}
\providecommand{\href}[2]{#2}

\end{document}